\documentclass[preprint]{elsarticle}
\usepackage{amsfonts}
\usepackage{amsmath}
\usepackage{amssymb}
\usepackage[all]{xy}
\usepackage{tikz}
\usetikzlibrary{trees}
\usetikzlibrary{arrows,shapes,snakes,automata,backgrounds,petri}

\newtheorem{thm}{Theorem}
\newtheorem{lem}[thm]{Lemma}
\newtheorem{cor}[thm]{Corollary}
\newtheorem{rem}{Remark}
\newproof{pf}{Proof}
\newproof{pot1}{Proof of Theorem \ref{mainthm}}

\begin{document}
\begin{frontmatter}

\title{A graph-theoretic condition for irreducibility of a set of cone preserving matrices}

\author[ref1]{Murad Banaji\corref{cor1}}
\author[ref1]{Andrew Burbanks}

\address[ref1]{Department of Mathematics, University of Portsmouth, Lion Gate Building, Lion Terrace, Portsmouth, Hampshire PO1 3HF, UK.}

\cortext[cor1]{Corresponding author: murad.banaji@port.ac.uk. }

\begin{abstract}
Given a closed, convex and pointed cone $K$ in $\mathbb{R}^n$, we present a result which infers $K$-irreducibility of sets of $K$-quasipositive matrices from strong connectedness of certain bipartite digraphs. The matrix-sets are defined via products, and the main result is relevant to applications in biology and chemistry. Several examples are presented. 
\end{abstract}

\begin{keyword}
partial order\sep convex cone\sep irreducibility\sep strong monotonicity

\MSC 15B48 \sep 15B35
\end{keyword}

\end{frontmatter}

\section{Introduction}
A digraph $G$ is strongly connected, or irreducible, if given any vertices $u$ and $v$, there exists a (directed) path from $u$ to $v$ in $G$. It is well known that a digraph is strongly connected if and only if its adjacency matrix is irreducible \cite{berman}. Here, given a cone $K$, we present a result which infers $K$-irreducibility of sets of $K$-quasipositive matrices from strong connectedness of associated bipartite digraphs. Graph-theoretic approaches to $K$-irreducibility of sets of $K$-positive matrices have been described in earlier work \cite{barkertam1981, barkertam1992, kunzesiegelpositivity}. These approaches are somewhat different in structure and philosophy to that described here. We will comment further on this in the concluding section. 

We will be interested in closed, convex cones in $\mathbb{R}^n$ which are additionally pointed (i.e., if $y \in K$ and $y \not = 0$, then $-y \not \in K$). Closed, convex and pointed cones will be abbreviated as {\bf CCP} cones. We do not assume the cones are solid (i.e., have nonempty interior in $\mathbb{R}^n$) -- however if a CCP cone is, additionally, solid, then it will be termed a {\bf proper} cone. For basic definitions and results on cones in $\mathbb{R}^n$ the reader is referred to \cite{berman,barker1973}. Given a CCP cone $K$, a {\bf face} $F \subseteq K$ will mean a closed face, namely $F$ is again a CCP cone, and moreover $x \in F$, $y \in K$, $x-y \in K$ together imply that $y \in F$. Faces other than $\{0\}$ and $K$ will be termed {\bf nontrivial}.

Let $K \subseteq\mathbb{R}^n$ be a CCP cone. Consider a real $n \times n$ matrix $M$. Recall that $M$ is {\bf $K$-positive} if $MK \subseteq K$. Defining $\mathbb{R}^n_{\geq 0}$ to be the (closed) nonnegative orthant in $\mathbb{R}^n$, a nonnegative matrix is then $\mathbb{R}^n_{\geq 0}$-positive. We will refer to $M$ as {\bf $K$-quasipositive} if there exists an $\alpha \in \mathbb{R}$ such that $M + \alpha I$ is $K$-positive. $\mathbb{R}^n_{\geq 0}$-quasipositive matrices -- generally referred to simply as quasipositive, or Metzler -- are just those with nonnegative off-diagonal entries. We define $M$ to be {\bf $K$-reducible} if there exists a nontrivial face $F$ of $K$ such that $M$ leaves $\mathrm{span}\,F$ invariant. This is a slight generalisation of the original definition of $K$-reducibility for $K$-positive matrices \cite{vandergraft1968} in order to allow us to apply the term to matrices which are not necessarily $K$-positive. A matrix which is not $K$-reducible is {\bf $K$-irreducible}. Note that an irreducible matrix could be termed $\mathbb{R}^n_{\geq 0}$-irreducible in this terminology. Alternatively any other orthant in $\mathbb{R}^n$ could be chosen as $K$.

\begin{rem}
\label{remposquasipos}
Given a CCP cone $K \subseteq\mathbb{R}^n$, an $n \times n$ matrix $M$ is $K$-irreducible if and only if $M + \alpha I$ is $K$-irreducible for each $\alpha \in \mathbb{R}$. In one direction we choose $\alpha = 0$. The other direction follows because given any face $F$ of $K$, $\mathrm{span}\,F$ is a vector subspace of $\mathbb{R}^n$. 
\end{rem}

{\bf Motivation from dynamical systems.} Motivation for examining $K$-irreducibility of a set of $K$-quasipositive matrices comes from the theory of monotone dynamical systems \cite{halsmith,hirschsmith}. Convex, pointed cones define partial orders in a natural way. Given a proper cone $K$ and a $C^1$ vector field $f:\mathbb{R}^n\to \mathbb{R}^n$, if the Jacobian matrix $Df(x)$ is $K$-quasipositive and $K$-irreducible at each $x \in \mathbb{R}^n$, then the associated local flow is strongly monotone with respect to the partial order defined by $K$. This result and a variety of technical modifications provide useful conditions which can be used to deduce the asymptotic behaviour of dynamical systems. 

\begin{rem}
Although the main results on monotone dynamical systems require the order cone to be solid, a CCP cone $K$ which fails to be solid is as useful as a proper one when the local flow or semiflow leaves cosets of $\mathrm{span}\,K$ invariant. Trivially, $K$ has nonempty relative interior in $\mathrm{span}\,K$, and attention can be restricted to cosets of $\mathrm{span}\,K$. This situation arises frequently in applications to biology and chemistry.
\end{rem}

\section{Some background material}

{\bf Sets of matrices.} It is convenient to use some notions from qualitative matrix theory. Let $M$ be a real matrix.
\begin{enumerate}
\item $\mathcal{Q}(M)$, the qualitative class of $M$ \cite{brualdi} is the set of all matrices with the same dimensions and sign pattern as $M$, namely if $N \in \mathcal{Q}(M)$, then $M_{ij} > 0 \Rightarrow N_{ij} > 0$, $M_{ij} < 0 \Rightarrow N_{ij} < 0$ and $M_{ij} = 0 \Rightarrow N_{ij} = 0$. 
\item $\mathcal{Q}_0(M)$ will stand for the closure of $\mathcal{Q}(M)$, namely if $N \in \mathcal{Q}_0(M)$, then $M_{ij} > 0 \Rightarrow N_{ij} \geq 0$, $M_{ij} < 0 \Rightarrow N_{ij} \leq 0$ and $M_{ij} = 0 \Rightarrow N_{ij} = 0$. Note that $\mathcal{Q}(M)$ and $\mathcal{Q}_0(M)$ are convex sets of matrices and $\mathcal{Q}(M) = \mathrm{relint}\,{\mathcal{Q}_0(M)}$, the relative interior of $\mathcal{Q}_0(M)$.
\item $\mathcal{Q}_1(M)$ will be the set of all matrices $N$ with the same dimensions as $M$ and satisfying $M_{ij}N_{ij} \geq 0$.
\end{enumerate}
Clearly, $\mathcal{Q}(M) \subseteq \mathcal{Q}_0(M) \subseteq \mathcal{Q}_1(M)$. 

\begin{rem}
\label{remorthant1}
Suppose $K$ is an orthant in $\mathbb{R}^n$ and $M$ is an $n \times n$ $K$-quasipositive matrix. It can easily be shown that each matrix in $\mathcal{Q}_0(M)$ is $K$-quasipositive. Since in this case $K$-irreducibility is simply irreducibility, if $M$ is $K$-irreducible then each matrix in $\mathcal{Q}(M)$ is $K$-irreducible. However, the same is not true for $\mathcal{Q}_0(M)$ which, after all, contains the zero matrix. 
\end{rem}

{\bf Notation for matrices.} Given any matrix $M$, we refer to the $k$th column of $M$ as $M_k$ and the $k$th row of $M$ as $M^k$. We define the new matrix $M^{(k)}$ by $M^{(k)}_{ij} = M_{ij}$ if $i = k$ and $M^{(k)}_{ij} = 0$ otherwise: i.e., $M^{(k)}$ is derived from $M$ by replacing all entries not in the $k$th row with zeros. Similarly $M_{(k)}$ is derived from $M$ by replacing all entries not in the $k$th column with zeros. A set of matrices $\mathbf{M}$ will be termed {\bf row-complete} if $M \in \mathbf{M} \Rightarrow M^{(k)} \in \mathbf{M}$ for each $k$, {\bf column-complete} if $M \in \mathbf{M} \Rightarrow M_{(k)} \in \mathbf{M}$ for each $k$, and {\bf complete} if $M \in \mathbf{M} \Rightarrow M^{(k)}, M_{(k)} \in \mathbf{M}$ for each $k$. Clearly, given any matrix $M$, $\mathcal{Q}_0(M)$ is complete; but smaller sets can be complete. For example, given some matrix $N$, 
\[
\mathbf{M} = \{D_1ND_2\,:\,D_1, D_2\,\,\mbox{are nonnegative diagonal matrices}\}
\]
is complete.

\begin{rem}
\label{remclasses}
Given a CCP cone $K$, a complete set of $K$-quasipositive matrices $\mathbf{M}$ and any $M \in \mathbf{M}$, $\mathcal{Q}_0(M)$ must consist of $K$-quasipositive matrices. This follows because (i) by completeness each matrix $M^{(i)}_{(j)}$ belongs to $\mathbf{M}$ and hence is $K$-quasipositive, and (ii) any finite nonnegative combination of $K$-quasipositive matrices is $K$-quasipositive. 
\end{rem}

{\bf Digraphs associated with square matrices.} Given an $n \times n$ matrix $M$, let $G_M$ be the associated digraph on $n$ vertices $u_1, \ldots, u_n$ defined in the usual way: the arc $u_iu_j$ exists in $G_M$ iff $M_{ij} \not = 0$. 

\begin{rem}
\label{remorthant2}
Following on from Remark~\ref{remorthant1}, if $K$ is an orthant in $\mathbb{R}^n$, then $K$-irreducibility of an $n \times n$ matrix $M$ is equivalent to strong connectedness of $G_{M}$.
\end{rem}

{\bf Bipartite digraphs associated with matrix-pairs.} Given an $n \times m$ matrix $A$ and an $m \times n$ matrix $B$, define a bipartite digraph $G_{A, B}$ on $n + m$ vertices as follows: associate a set of $n$ vertices $u_1, \ldots, u_n$ with the rows of $A$ (we will refer to these as the ``row vertices'' of $G_{A, B}$); associate another $m$ vertices $v_1, \ldots, v_m$ with the columns of $A$ (we will refer to these as the ``column vertices'' of $G_{A, B}$); add the arc $u_iv_j$ iff $A_{ij} \not = 0$; add the arc $v_ju_i$ iff $B_{ji} \not = 0$. 

\begin{rem}
This is a specialisation of the general construction of block-circulant digraphs from sets of appropriately dimensioned matrices in \cite{banajirutherford1}. If $A$ and $B$ are $(0,1)$ matrices (with $AB$ a square matrix), then the adjacency matrix of $G_{A, B}$ is simply
\[
\left(\begin{array}{cc}0 & A\\B & 0\end{array}\right).
\]
The bipartite digraph here is also closely related to the so-called DSR graph, presented in \cite{banajicraciun2} and used to make claims about dynamical systems arising in biology and chemistry.
\end{rem}

{\bf The context of the main result.} Fundamental early results on convergence in monotone dynamical systems \cite{hirsch85} apply to systems with Jacobian matrices which are quasipositive and irreducible (in our terminology $\mathbb{R}^n_{\geq 0}$-quasipositive and $\mathbb{R}^n_{\geq 0}$-irreducible). Generalising from $\mathbb{R}^n_{\geq 0}$ to all orthants is straightforward: where $K$ is an orthant, there is a simple graph-theoretic test \cite{halsmith} to decide $K$-quasipositivity of a given $n \times n$ matrix $M$. By Remarks~\ref{remorthant1}~and~\ref{remorthant2}, $K$-quasipositivity extends to all of $\mathcal{Q}_0(M)$, and $K$-irreducibility of some $M' \in \mathcal{Q}_0(M)$ is equivalent to strong connectedness of $G_{M'}$.

We are interested in how this situation generalises when $K$ is not necessarily an orthant, not necessarily simplicial, and in fact not necessarily finitely generated. In general, given a $K$-quasipositive matrix $M$, we can rarely expect all matrices in $\mathcal{Q}_0(M)$ to be $K$-quasipositive. However the following situation is not uncommon: there are matrices $A$ and $\tilde B$ such that $AB$ is $K$-quasipositive for each $B\in \mathcal{Q}_0(\tilde B)$. The practical relevance is to applications in biology and chemistry where Jacobian matrices often have a constant initial factor, but a second factor with variable entries whose signs are, however, known. A number of particular examples were presented in \cite{banajidynsys}.  

Given a set of $K$-quasipositive matrices of the form $\{AB\,:\, B \in \mathbf{B}\}$, we would hope that there is a natural graph-theoretic test to decide which members of this set are also $K$-irreducible. Theorem~\ref{mainthm} below provides precisely such a condition: provided $\mathbf{B}$ is complete and the initial factor $A$ satisfies a mild genericity condition, $K$-irreducibility of $AB$ follows from strong connectedness of the bipartite digraph $G_{A, B}$. In the special case where $K$ is the nonnegative orthant, $A$ is the identity matrix, and $\mathbf{B}$ is the set of nonnegative matrices, the results reduce to well known ones.
\begin{rem}
Our motivation for considering complete sets of matrices will be as follows. Consider a set of $K$-quasipositive matrices $\mathcal{M} = \{AB\,:\, B \in \mathbf{B}\}$, where $A$ is $n \times m$, $\mathbf{B}$ consists of $m \times n$ matrices, and $\mathbf{B}$ is complete. Clearly
\[
AB = A\sum_k B^{(k)} = \sum_k AB^{(k)} = \sum_k A_kB^k.
\]
So, by row-completeness, any matrix in $\mathcal{M}$ can be written as a sum of rank $1$ $K$-quasipositive matrices in $\mathcal{M}$. On the other hand, suppose, for some $v \in \mathbb{R}^n$, some $i \in \{1, \ldots, m\}$, and all $B \in \mathbf{B}$, that $(Bv)_i \geq 0$ (resp. $(Bv)_i \leq 0$). Then, by column-completeness, for each $B$ and each $k \in \{1, \ldots, n\}$, $(B_{(k)}v)_i =B_{ik}v_k \geq 0$ (resp. $B_{ik}v_k \leq 0$), i.e., $B^i \in \mathcal{Q}_1(v^{\mathrm{T}})$ (resp. $B^i \in \mathcal{Q}_1(-v^{\mathrm{T}})$).
\end{rem}

\section{The main result}

From now on $K \subseteq\mathbb{R}^n$ will be a CCP cone in $\mathbb{R}^n$, $A$ an $n \times m$ matrix and $\mathbf{B}$ a complete set of $m \times n$ matrices. For any $B \in \mathbf{B}$, $AB$ is an $n \times n$ matrix. The main result of this paper is:

\begin{thm}
\label{mainthm}
Assume that $\mathrm{Im}\,A \not \subseteq \mathrm{span}\,F$ for any nontrivial face $F$ of $K$. Suppose that for each $B \in \mathbf{B}$, $AB$ is $K$-quasipositive. Then whenever $G_{A, B}$ is strongly connected, $AB$ is also $K$-irreducible. 
\end{thm}
An immediate corollary is:
\begin{cor}
\label{maincor}
Assume that $\mathrm{Im}\,A \not \subseteq \mathrm{span}\,F$ for any nontrivial face $F$ of $K$. Suppose that for each $B \in \mathbf{B}$, $AB$ is $K$-positive. Then whenever $G_{A, B}$ is strongly connected, $AB$ is also $K$-irreducible. 
\end{cor}
\begin{pf}
$K$-positivity of $AB$ implies $K$-quasipositivity of $AB$. The result now follows from Theorem~\ref{mainthm}. \qquad \qed
\end{pf}

\begin{rem}
Note that if $\mathrm{Im}\,A \subseteq \mathrm{span}\,F$ for some nontrivial face $F$ of $K$, then trivially $AB$ is $K$-reducible. To see that the assumption that $\mathrm{Im}\,A \not \subseteq \mathrm{span}\,F$ is in general necessary in Theorem~\ref{mainthm}, consider the matrices
\[
\Lambda = \left(\begin{array}{cc}1&2\\2&1\end{array}\right), \quad A = \left(\begin{array}{cc}1&1\\2&2\end{array}\right), \quad B = \left(\begin{array}{cc}a&b\\c&d\end{array}\right)
\]
where $a,b,c,d \geq 0$. Let $K = \{\Lambda z\,:\, z \in \mathbb{R}^2_{\geq 0}\}$. Then
\[
AB\Lambda_1 = (a+c + 2(b+d))\Lambda_1, \quad AB\Lambda_2 = (2(a+c) + b+d)\Lambda_1
\]
which are both clearly in $K$ for any $a,b,c,d \geq 0$. So $AB$ is $K$-positive. On the other hand $F = \{r\Lambda_1\,:\,r \geq 0\}$ satisfies $(AB)F \subseteq F$ for all $B$, so $AB$ is $K$-reducible for all $a, b, c, d \geq 0$. However, for $a, b, c, d > 0$, $G_{A, B}$ is a complete bipartite digraph which is obviously strongly connected. 
\end{rem}

\section{Proofs}

We need some preliminary lemmas in order to prove Theorem~\ref{mainthm}. The following is proved as Lemma~4.4 in \cite{banajidynsys}:
\begin{lem}
\label{plusminus}
Let $F$ be a face of $K$, $v_1, v_2 \in F$, and $w \in \mathbb{R}^n$. If there exist $\alpha, \beta > 0$ such that $v_1 + \alpha w \in K$ and $v_2 - \beta w \in K$, then $w \in \mathrm{span}(F)$. 
\end{lem}
\begin{pf}
Define $y_1 \equiv v_1 + \alpha w$ and $y_2 \equiv v_2 - \beta w$. Then $y_3 \equiv y_1 + (\alpha/\beta)y_2 = v_1 + (\alpha/\beta)v_2 \in F$. Since $y_3 \in F$, $y_3 - y_1 = (\alpha/\beta)y_2 \in K$, and $y_1 \in K$, by the definition of a face, $y_1 \in F$. So $w = (y_1 - v_1)/\alpha \in \mathrm{span}(F)$.
\qquad \qed \end{pf}
{\bf Extremals.} A one dimensional face of $K$ will be termed an extremal ray or an extremal for short, and any nonzero vector in an extremal ray will be an extremal vector of $K$. 

\begin{lem}
\label{qualclass}
Suppose that for each $B \in \mathbf{B}$, $AB$ is $K$-quasipositive. Let $v$ be an extremal vector of $K$. Then for each $j$ either $A_j = rv$ for some real number $r$ or $B^j \in \mathcal{Q}_1(v^\mathrm{T})$ for all $B \in \mathbf{B}$, or $B^j \in \mathcal{Q}_1(-v^\mathrm{T})$ for all $B \in \mathbf{B}$. 
\end{lem}
\begin{pf}
Suppose there exist $j$ and $P,Q \in \mathbf{B}$ such that ${P}^j \not \in \mathcal{Q}_1(-v^\mathrm{T})$ and ${Q}^j \not \in \mathcal{Q}_1(v^\mathrm{T})$, i.e., there exist $k, l$ such that $P_{jk}v_k \equiv t_1 > 0$,  $-Q_{jl}v_l \equiv t_2  > 0$. ($P=Q$ is possible.) Defining $\overline{B} = P_{(k)}$, $\underline{B} = Q_{(l)}$, note that since $\mathbf{B}$ is column-complete, $\overline{B},\underline{B} \in \mathbf{B}$. By construction, $\overline{B}^j v  = t_1 > 0$ and $-\underline{B}^j v = t_2 > 0$. Since $\mathbf{B}$ is row-complete, $\overline{B}^{(j)}, \underline{B}^{(j)} \in \mathbf{B}$ and so $A\overline{B}^{(j)}$, $A\underline{B}^{(j)}$ are $K$-quasipositive. Let $\alpha_1$ and $\alpha_2$ be such that $A\overline{B}^{(j)} v + \alpha_1v \in K$ and $A\underline{B}^{(j)} v + \alpha_2v \in K$ respectively. We can assume (w.l.o.g.) that $\alpha_1, \alpha_2 > 0$. Then:
\[
A\overline{B}^{(j)} v + \alpha_1v  = A_j\overline{B}^j v + \alpha_1v  = \alpha_1v + t_1A_j \in K
\]
and 
\[
A\underline{B}^{(j)} v + \alpha_2v = A_j\underline{B}^j v + \alpha_2v = \alpha_2v -t_2A_j \in K. 
\]
By Lemma~\ref{plusminus}, these two equations imply that $A_j = rv$ for some $r$. 
\qquad \qed
\end{pf}

\begin{lem}
\label{decomp}
Let $F$ be a nontrivial face of $K$, $\{J(k)\}$ be a finite set of $n \times n$ $K$-quasipositive matrices and $J = \sum J(k)$. If there exists $x \in F$ such that $Jx \in \mathrm{span}\,F$, then $J(k) x \in \mathrm{span}\,F$ for each $k$.
\end{lem} 
\begin{pf}
By $K$-quasipositivity of each $J(k)$ we can write
\[
J(k)x = p_k + q_k
\]
where $p_k \in (K\backslash F) \cup \{0\}$ and $q_k \in \mathrm{span}\,F$. Summing, we get
\[
Jx = p +q
\]
where $q  = \sum_kq_k \in \mathrm{span}\,F$ and $p = \sum_kp_k \in (K\backslash F) \cup \{0\}$. Now if $Jx \in \mathrm{span}\,F$ then $p = 0$. Since $p_k \in K$ and $K$ is pointed, this implies that $p_k = 0$ for each $k$. So $J(k) x = q_k\in \mathrm{span}\,F$ for each $k$, proving the lemma. \qquad \qed
\end{pf}

\begin{lem}
\label{lemtechnical}
Let $AB$ be $K$-quasipositive for each $B \in \mathbf{B}$ and let $F$ be a nontrivial face of $K$ spanned by (pairwise independent) extremal vectors $\{\Lambda_i\}_{i \in \mathcal{I}}$. Choose and fix some $B\in \mathbf{B}$. Then given any nonempty set $\mathcal{R} \subseteq \{1, \ldots, m\}$, either  (i) $A_k \in \mathrm{span}\,F$ for some $k \in \mathcal{R}$ or (ii) $B^k\Lambda_i = 0$ for each $k \in \mathcal{R}, i \in \mathcal{I}$, or (iii) there exists $\Lambda_i \in F$ such that $\sum_{k \in \mathcal{R}}A_kB^k\Lambda_i \not \in \mathrm{span}\,F$. 
\end{lem} 

\begin{pf}
For each $k$, recall that $B^{(k)} \in \mathbf{B}$, so $AB^{(k)} = A_kB^k$ is $K$-quasipositive. Suppose (iii) does not hold, i.e., $\sum_{k \in \mathcal{R}}A_kB^k\Lambda_i \in \mathrm{span}\,F$ for each $\Lambda_i \in F$. Applying Lemma~\ref{decomp} with $J(k) = A_kB^k$, we get for each $k \in \mathcal{R}$ and each $\Lambda_i \in F$ that $A_kB^k\Lambda_i \in \mathrm{span}\,F$. So for each fixed $k \in \mathcal{R}$, either $A_k \in \mathrm{span}\,F$ or $B^k\Lambda_i=0$ for all $i \in \mathcal{I}$. \qquad \qed \end{pf}

\begin{pot1}
We show that if $AB$ is $K$-reducible for some $B \in \mathbf{B}$, then $G_{A, B}$ cannot be strongly connected. Let $\{\Lambda_i\}$ be a set of pairwise independent extremal vectors generating $K$. Let $F$ be a nontrivial face of $K$ such that $\mathrm{span}\,F$ is left invariant by $AB$, and define $\mathcal{I}$ via $i \in \mathcal{I}\Leftrightarrow \Lambda_i \in F$. Let $\mathcal{R} \subsetneq \{1, \ldots, m\}$ be defined by $k \in \mathcal{R} \Leftrightarrow A_k \in \mathrm{span}\,F$. $\mathcal{R}$ may be empty, but by assumption cannot be all of $\{1, \ldots, m\}$ since $\mathrm{Im}\,A \not \subseteq \mathrm{span}\,F$. So $\mathcal{R}^c$, the complement of $\mathcal{R}$, is nonempty. 

Choose any $x \in F$. We have
\begin{equation}
\label{eq3terms}
AB x = \sum_k A_kB^kx= \sum_{k \in \mathcal{R}}A_kB^kx + \sum_{k \in \mathcal{R}^c}A_kB^kx\,.
\end{equation}
Clearly $\sum_{k \in \mathcal{R}}A_kB^kx \in \mathrm{span}\,F$. By assumption, $AB x \in \mathrm{span}\,F$,  and so $\sum_{k \in \mathcal{R}^c}A_kB^kx \in \mathrm{span}\,F$. Now since $x \in F$ was arbitrary and $A_k \not \in \mathrm{span}\,F$ for any $k \in \mathcal{R}^c$, by Lemma~\ref{lemtechnical} we must have $B^k\Lambda_i = 0$ for each $k \in \mathcal{R}^c, i \in \mathcal{I}$. But from Lemma~\ref{qualclass} we know that either (i) $A_k = r\Lambda_i$ for some scalar $r$ or (ii) $B^k \in \mathcal{Q}_1(\Lambda_i^\mathrm{T})$ or $B^k \in \mathcal{Q}_1(-\Lambda_i^\mathrm{T})$. Since $A_k \not \in \mathrm{span}\,F$ the first possibility is ruled out, and (ii) must hold. Consequently $B^k \Lambda_i = 0$ implies $B_{kl} \Lambda_{li} = 0$ for each $l$. 

The above is true for each $i \in \mathcal{I}$, $k \in \mathcal{R}^c$. Now there are two possibilities:
\begin{enumerate}
\item Suppose that $\mathcal{R}$ is empty. Then, for each $i \in \mathcal{I}$, and all $k, l$, $B_{kl} \Lambda_{li} = 0$. Fix some $i \in \mathcal{I}$ and some $l$ such that $\Lambda_{li} \not = 0$; then $B_{kl} = 0$ for all $k$ (the $l$th column of $B$ is zero). By the definition of $G_{A, B}$, this means that there are no arcs incident into the row vertex $u_l$. Thus $G_{A, B}$ is not strongly connected.
\item Suppose that $\mathcal{R}$ is nonempty. Given $k' \in \mathcal{R}$ we can write $A_{k'} = \sum_{i \in \mathcal{I}}q_i\Lambda_i$ for some constants $q_i$, so for any $k \in \mathcal{R}^c$,
\[
B^kA_{k'} = \sum_{i \in \mathcal{I}}q_iB^k\Lambda_{i} = 0.
\]
Moreover, since for each $i \in \mathcal{I}$ and each $l$, $B_{kl} \Lambda_{li} = 0$, we get $B_{kl} A_{lk'} = 0$. In terms of $G_{A, B}$, this means that there is no directed path of length $2$ from any column vertex $v_k$ with $k \in \mathcal{R}^c$ to a column vertex $v_{k'}$ with $k' \in\mathcal{R}$. Thus there is no directed path (of any length) of the form $v_k\cdots v_{k'}$ with $k \in \mathcal{R}^c$, $k' \in \mathcal{R}$, and $G_{A, B}$ is not strongly connected. 
\end{enumerate}
\qquad\qed
\end{pot1}

\section{Examples}
The examples below illustrate application of Theorem~\ref{mainthm} and Corollary~\ref{maincor}.\\

{\em Example 1.} Consider the special case where Corollary~\ref{maincor} is applied with $A = I$, $\mathbf{B}$ the nonnegative matrices and $K = \mathbb{R}^n_{\geq 0}$. Since $\mathrm{Im}\,I = \mathbb{R}^n$, clearly $\mathrm{Im}\,I \not \subseteq \mathrm{span}\,F$ for any nontrivial face $F$ of $\mathbb{R}^n_{\geq 0}$. It is also immediate that $\mathbf{B}$ is complete and for each $B \in \mathbf{B}$, the $n \times n$ matrix $IB$ is $K$-positive. Now we show that $G_{I, B}$ is strongly connected if and only if $G_B$ is strongly connected. (i) Suppose $G_B$ is strongly connected. An arc from vertex $i$ to vertex $j$ in $G_B$ implies that $B_{ij} \not = 0$. But $B_{ij} = I_{ii}B_{ij}$, and since $I_{ii} = 1$ this implies that there exists a path $u_iv_iu_j$ in $G_{I, B}$. Thus a path from vertex $i$ to vertex $j$ in $G_B$ implies a path from vertex $u_i$ to vertex $u_j$ in $G_{I, B}$. Thus strong connectedness of $G_B$ implies a path between any two row vertices $u_i$ and $u_j$ of $G_{I, B}$. On the other hand since all arcs $u_iv_i$ exist in $G_{I, B}$, the path $u_i\cdots u_j$ implies the existence of paths $u_i\cdots v_j$, $v_i \cdots u_j$ and $v_i \cdots v_j$. (ii) Suppose $G_{I, B}$ is strongly connected. The path  $v_i \cdots v_j$ in $G_{I, B}$ immediately implies a path from vertex $i$ to vertex $j$ in $G_B$. Thus we recover from Corollary~\ref{maincor} the fact that for a nonnegative matrix $B$, strong connectedness of $G_B$ implies irreducibility of $B$.\\

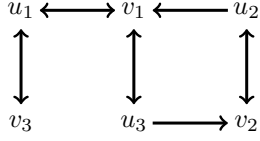
\begin{figure}[h]
\begin{minipage}{\textwidth}
\begin{center}
\begin{tikzpicture}[domain=0:4,scale=0.5]

\node at (1,4) {$u_1$};
\node at (4,4) {$v_1$};
\node at (7,4) {$u_2$};
\node at (1,1) {$v_3$};
\node at (4,1) {$u_3$};
\node at (7,1) {$v_2$};

\draw[<->, line width=0.04cm] (1.5, 4) -- (3.5,4);
\draw[<-, line width=0.04cm] (4.5, 4) -- (6.5,4);
\draw[->, line width=0.04cm] (4.5, 1) -- (6.5,1);

\draw[<->, line width=0.04cm] (1, 3.5) -- (1,1.5);
\draw[<->, line width=0.04cm] (4, 3.5) -- (4,1.5);
\draw[<->, line width=0.04cm] (7, 3.5) -- (7,1.5);

\end{tikzpicture}
\end{center}
\end{minipage}
\caption{\label{DSR1} $G_{A, B}$ for the system in Example 2 with $B \in \mathrm{relint}\,\mathbf{B}$. By inspection the digraph is strongly connected.}
\end{figure}

{\em Example 2.} Let
\[
A = \Lambda = 
\left(
\begin{array}{rrr}
1 & 0 & -1 \\
1 & -2 & 0 \\
-1 & 1 & 0
\end{array}
\right) \quad \mbox{and} \quad 
\tilde B = 
\left(
\begin{array}{rrr}
-1 & 0 & 1 \\
0 & 1 & 0 \\
1 & 0 & 0
\end{array}
\right)
\]
Let $\mathbf{B} = \mathcal{Q}_0(\tilde B)$ (so $\mathbf{B}$ is complete) and define 
\[
K = \{\Lambda z\,:\, z \in \mathbb{R}^3_{\geq 0}\}.
\]
$K$ is a CCP cone in $\mathbb{R}^3$. (It is easy to show that any nonsingular $n \times n$ matrix defines a proper simplicial cone in $\mathbb{R}^n$ in this way.) Since $A = \Lambda$ it is immediate that $\mathrm{Im}\,A$ does not lie in the span of any nontrivial face of $K$. Consider any $B \in \mathbf{B}$, i.e., any matrix of the form
\[
B = 
\left(
\begin{array}{rrr}
-a & 0 & b \\
0 & c & 0 \\
d & 0 & 0
\end{array}
\right)
\] 
where $a,b,c,d \geq 0$. Since $A = \Lambda$,
\[
AB\Lambda + (a+b+2c+d)\Lambda = \Lambda(B\Lambda + (a+b+2c+d) I)
\]
and it can be checked that $B\Lambda + (a+b+2c+d) I$ is nonnegative. Thus $AB$ is $K$-quasipositive for all $B \in \mathbf{B}$. On the other hand $G_{A, B}$ is illustrated in Figure~\ref{DSR1} for any $B \in \mathrm{relint}\,\mathbf{B}$ (i.e. $B \in \mathcal{Q}(\tilde B)$) and can be seen to be strongly connected. So $AB$ is $K$-irreducible for any $B \in \mathrm{relint}\,\mathbf{B}$. \\

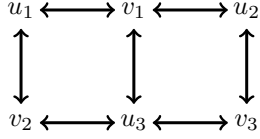
\begin{figure}[h]
\begin{minipage}{\textwidth}
\begin{center}
\begin{tikzpicture}[domain=0:4,scale=0.5]

\node at (1,1) {$v_2$};
\node at (4,1) {$u_3$};
\node at (7,1) {$v_3$};
\node at (1,4) {$u_1$};
\node at (4,4) {$v_1$};
\node at (7,4) {$u_2$};

\draw[<->, line width=0.04cm] (1.5, 4) -- (3.5,4);
\draw[<->, line width=0.04cm] (1.5, 1) -- (3.5,1);
\draw[<->, line width=0.04cm] (4.5, 4) -- (6.5,4);
\draw[<->, line width=0.04cm] (4.5, 1) -- (6.5,1);

\draw[<->, line width=0.04cm] (1, 3.5) -- (1,1.5);
\draw[<->, line width=0.04cm] (4, 3.5) -- (4,1.5);
\draw[<->, line width=0.04cm] (7, 3.5) -- (7,1.5);

\end{tikzpicture}
\end{center}
\end{minipage}
\caption{\label{DSR2} $G_{A, B}$ for the system in Example 3 with $B \in \mathrm{relint}\,\mathbf{B}$. By inspection the digraph is strongly connected.}
\end{figure}

{\em Example 3.} The following is an example with a cone which is not solid. Let
\[
A = 
\left(
\begin{array}{rrr}
-1 & -1 & 0 \\
-1 & 0 & -1 \\
2 & 1 & 1
\end{array}
\right), \quad \Lambda = 
\left(
\begin{array}{rr}
1 & 0 \\
0 & -1 \\
-1 & 1
\end{array}
\right) \quad \mbox{and} \quad 
\tilde B = 
\left(
\begin{array}{rrr}
1 & 1 & -1 \\
1 & 0 & -1 \\
0 & 1 & -1
\end{array}
\right)
\]
Let $\mathbf{B} = \mathcal{Q}_0(\tilde B)$ and define 
\[
K = \{\Lambda z\,:\, z \in \mathbb{R}^2_{\geq 0}\}.
\]
$K$ is a CCP cone of dimension $2$ in $\mathbb{R}^3$. Clearly $\mathrm{Im}\,A$ does not lie in the span of any nontrivial face (i.e., any extremal) of $K$. Defining 
\[
T = 
\left(
\begin{array}{rrr}
-1 & -1 & 0 \\
1 & 0 & 1
\end{array}
\right),
\]
note that $A = \Lambda T$. Consider any $B \in \mathbf{B}$, i.e., any matrix of the form
\[
B = 
\left(
\begin{array}{rrr}
a & b & -c \\
d & 0 & -e \\
0 & f & -g
\end{array}
\right)
\] 
where $a,b,c,d,e,f,g \geq 0$. Then
\[
AB\Lambda + (a+b+c+d+e+f+g)\Lambda = \Lambda(TB\Lambda + (a+b+c+d+e+f+g) I)
\]
where $I$ is the $2 \times 2$ identity matrix. It can be checked that $TB\Lambda + (a+b+c+d+e+f+g) I$ is nonnegative. Thus $AB$ is $K$-quasipositive for all $B \in \mathbf{B}$. $G_{A, B}$ is illustrated in Figure~\ref{DSR2} for any $B \in \mathrm{relint}\,\mathbf{B}$ (i.e., $B \in \mathcal{Q}(\tilde B)$) and can be seen to be strongly connected. So $AB$ is $K$-irreducible for any $B \in \mathrm{relint}\,\mathbf{B}$.\\

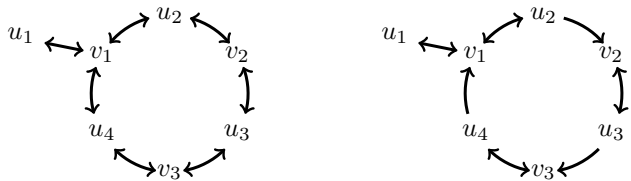
\begin{figure}[h]
\begin{center}
\begin{minipage}{0.4\textwidth}
\begin{tikzpicture}[domain=0:12,scale=0.52]

\path (90: 2cm) coordinate (A2);
\path (30: 2cm) coordinate (R2);
\path (-30: 2cm) coordinate (A3);
\path (270: 2cm) coordinate (R3);
\path (210: 2cm) coordinate (A1);
\path (150: 2cm) coordinate (R1);

\path (-3.8, 1.5) coordinate (B1);

\node at (A1) {$u_4$};
\node at (A2) {$u_2$};
\node at (A3) {$u_3$};

\node at (R1) {$v_1$};
\node at (R2) {$v_2$};
\node at (R3) {$v_3$};

\node at (B1) {$u_1$};

\path (90-15: 2cm) coordinate (A2end);
\path (-30-15: 2cm) coordinate (A3end);
\path (210-15: 2cm) coordinate (A1end);

\path (90+15: 2cm) coordinate (A2end1);
\path (-30+15: 2cm) coordinate (A3end1);
\path (210+15: 2cm) coordinate (A1end1);

\draw[<->, line width=0.04cm] (A2end)  arc (90-15:90-50:2cm);
\draw[<->, line width=0.04cm] (A3end)  arc (-30-15:-30-50:2cm);
\draw[<->, line width=0.04cm] (A1end)  arc (210-15:210-50:2cm);

\draw[<->, line width=0.04cm] (A2end1)  arc (90+15:90+50:2cm);
\draw[<->, line width=0.04cm] (A3end1)  arc (-30+15:-30+50:2cm);
\draw[<->, line width=0.04cm] (A1end1)  arc (210+15:210+50:2cm);

\draw[<->, line width=0.04cm] (-2.2, 1.1) -- (-3.2, 1.3);

\end{tikzpicture}
\end{minipage}
\begin{minipage}{0.4\textwidth}
\begin{tikzpicture}[domain=0:12,scale=0.52]


\path (90: 2cm) coordinate (A2);
\path (30: 2cm) coordinate (R2);
\path (-30: 2cm) coordinate (A3);
\path (270: 2cm) coordinate (R3);
\path (210: 2cm) coordinate (A1);
\path (150: 2cm) coordinate (R1);

\path (-3.8, 1.5) coordinate (B1);

\node at (A1) {$u_4$};
\node at (A2) {$u_2$};
\node at (A3) {$u_3$};

\node at (R1) {$v_1$};
\node at (R2) {$v_2$};
\node at (R3) {$v_3$};

\node at (B1) {$u_1$};

\path (90-15: 2cm) coordinate (A2end);
\path (-30-15: 2cm) coordinate (A3end);
\path (210-15: 2cm) coordinate (A1end);

\path (90+15: 2cm) coordinate (A2end1);
\path (-30+15: 2cm) coordinate (A3end1);
\path (210+15: 2cm) coordinate (A1end1);

\draw[->, line width=0.04cm] (A2end)  arc (90-15:90-50:2cm);
\draw[->, line width=0.04cm] (A3end)  arc (-30-15:-30-50:2cm);
\draw[->, line width=0.04cm] (A1end)  arc (210-15:210-50:2cm);

\draw[<->, line width=0.04cm] (A2end1)  arc (90+15:90+50:2cm);
\draw[<->, line width=0.04cm] (A3end1)  arc (-30+15:-30+50:2cm);
\draw[<->, line width=0.04cm] (A1end1)  arc (210+15:210+50:2cm);

\draw[<->, line width=0.04cm] (-2.2, 1.1) -- (-3.2, 1.3);

\end{tikzpicture}
\end{minipage}
\end{center}

\caption{\label{DSR} {\em Left.} The digraph $G_{A, B}$ where $A$ is as shown in (\ref{exampleeqn}) and $B$ is any matrix in $\mathcal{Q}(-A^\mathrm{T})$. $u_i$ is the row vertex corresponding to row $i$ in $A$, while $v_i$ corresponds to column $i$. {\em Right.} Setting $B_{22} = B_{33} = B_{14} = 0$ removes the arcs $v_2u_2$, $v_3u_3$, $v_1u_4$ from $G_{A, B}$, but still gives a strongly connected digraph.}
\end{figure}

{\em Example 4.} As a final, nontrivial, example, let 
\begin{equation}
\label{exampleeqn}
A = \left(\begin{array}{rrr}-1&0&0\\1&-1&0\\0&1&-1\\1&0&1\end{array}\right), \quad \Lambda =  \left(\begin{array}{rrrrrrrr}\,\,0&\,\,0&\,\,0&0&1&1&1&1\\0&0&1&1&-1&-1&0&0\\1&0&0&-1&1&0&0&-1\\0&1&0&1&-1&0&-1&0\end{array}\right),
\end{equation}
and $\mathbf{B} = \mathcal{Q}_0(-A^\mathrm{T})$. Define 
\[
K = \{\Lambda z\,:\, z \in \mathbb{R}^8_{\geq 0}\}.
\]
Various facts can be confirmed either theoretically or via computation:
\begin{enumerate}
\item $K$ is a proper cone in $\mathbb{R}^4$.
\item $\mathrm{Im}\,A \not \subseteq \mathrm{span}\,F$ for any nontrivial face $F$ of $K$.
\item For each $B \in \mathbf{B}$, $AB$ is $K$-quasipositive.
\end{enumerate}
Some insight into the structure of $K$ and the proof of these facts is provided in the Appendix. It now follows from Theorem~\ref{mainthm} that whenever $G_{A, B}$ is strongly connected, $AB$ (and hence $AB + \alpha I$ for each $\alpha \in \mathbb{R}$) is also $K$-irreducible. For example, it can easily be checked that for $B \in \mathrm{relint}\,\mathbf{B}$ (namely $B \in \mathcal{Q}(-A^\mathrm{T})$), $G_{A, B}$ is strongly connected and so $AB$ is $K$-irreducible. The condition that $B \in \mathrm{relint}\,\mathbf{B}$ can be relaxed while maintaining $K$-irreducibility. The digraphs $G_{A, B}$ for two choices of $B \in \mathbf{B}$ are illustrated in Figure~\ref{DSR}.

Confirming $K$-irreducibility for choices of $B$ without the aid of Theorem~\ref{mainthm} is possible but tedious, requiring computation of the action of $AB$ on each of the 26 nontrivial faces of $K$. 

\section{Concluding remarks}

Rather different graph-theoretic approaches to questions of irreducibility of matrices are taken in \cite{barkertam1981, barkertam1992, kunzesiegelpositivity}. In \cite{kunzesiegelpositivity}, for example, polyhedral cones with $n_K$ extremals were considered, and digraphs on $n_K$ vertices constructed. Results were presented deriving $K$-irreducibility of matrices from $K$-quasipositivity of these matrices and strong connectedness of the digraphs. The construction relies, however on knowledge of the matrix action on each extremal of $K$.

In our approach described here, $K$ is not necessarily polyhedral and no knowledge of the facial structure or particular action of matrices on extremals of $K$ is required. In compensation, however, we assume that a set of matrices with a particular structure (completeness) are {\em all} $K$-quasipositive, and $K$-quasipositivity of this entire set is essential for the proofs. This stronger assumption about $K$-quasipositivity allows weaker assumptions about the structure of $K$ and the action of the matrices on faces of $K$. Thus both the construction of the digraph here and the assumptions are somewhat different from earlier work in this area.

\section*{Acknowledgements}

The work of M. Banaji was supported by EPSRC grant EP/J008826/1 ``Stability and order preservation in chemical reaction networks''.

\appendix

\section{Some details connected with Example 4}
{\bf $K$ is a proper cone.} That $K$ is closed and convex is immediate from the definition. Note that $\Lambda P = I$ where $P$ is the nonnegative matrix
\[
P = \left(\begin{array}{cccc}0&0&1&0\\1&0&0&1\\0&1&0&0\\0&0&0&0\\0&0&0&0\\0&0&0&0\\1&0&0&0\\0&0&0&0\end{array}\right).
\]
Consequently $K$ contains $\mathbb{R}^4_{\geq 0}$ and $K$ is solid. Defining $p = (2,1,1,1)^\mathrm{T}$, $\Lambda^{\mathrm{T}}p$ is strictly positive, so $p^\mathrm{T}\Lambda z > 0$ for any nonnegative and nonzero $z$. Thus $p \in \mathrm{int}\,K^*$, where $K^*$ is the dual cone to $K$. Since $K^*$ has nonempty interior this implies that $K$ is pointed (if $K$ contains a nonzero $y \in \mathbb{R}^4$ such that $y, -y \in K$, we get the contradiction $p^\mathrm{T}y > 0$ and $p^\mathrm{T}(-y) > 0$). Putting together these observations, $K$ is a proper cone in $\mathbb{R}^4$. \\

{\bf The facial structure of $K$.} It can be checked that each $\Lambda_i$ spans a different extremal of $K$, namely, no $\Lambda_i$ can be constructed as a nonnegative combination of others. Further, the two dimensional faces of $K$ are spanned by pairs of $\Lambda_i$ for $i$ belonging to:
\[
\{1,2\}, \{1,3\}, \{2,4\},\{3,4\},\{1,5\},\{2,6\},\{3,7\},\{4,8\},\{5,6\},\{5,7\},\{6,8\},\{7,8\}.
\]
while the three dimensional faces of $K$ are spanned by sets of four $\Lambda_i$ for $i$ belonging to:
\[
\{1,2,3,4\}, \{1,2,5,6\},\{3,4,7,8\},\{1,3,5,7\},\{2,4,6,8\} \,\,\, \mbox{and} \,\,\, \{5,6,7,8\}.
\]
This facial structure of $K$ can be computed directly, or is deduced as follows. Let 
\[
\Gamma = \left(\begin{array}{rrr}0&0&1\\0&1&-1\\-1&-1&0\\1&0&-1\end{array}\right)\, \quad \mbox{and} \quad C = \left(\begin{array}{cccccccc}0&1&0&1&0&1&0&1\\0&0&1&1&0&0&1&1\\0&0&0&0&1&1&1&1\end{array}\right)\,.
\]
Note that the columns of $C$ are the vertices of a cube $\mathcal{C} \subseteq \mathbb{R}^3$ and the index sets above define the one and two dimensional faces of $\mathcal{C}$. Define $q = (0,0,1,0)^{\mathrm{T}}$ and the affine mapping $T:\mathbb{R}^3 \to \mathbb{R}^4$ by $T(v) = \Gamma v+q$. Observe that $\Lambda_i = T(C_i)$, that $\Gamma$ has trivial kernel, and that $q \not \in \mathrm{Im}\,\Gamma$. As noted above, $\Lambda_i$ are pairwise independent extremal vectors generating $K$. It is easy to confirm that as a consequence each $x \in K\backslash\{0\}$ can be written $x = rT(y)$, where $r$ and $y$ are uniquely defined.
Given any nonempty $\mathcal{I} \subseteq \{1, \ldots, 8\}$, define $F_{\mathcal{I}} \subseteq \mathcal{C}$ as the convex hull of $\{C_i\}_{i \in \mathcal{I}}$, and $F'_{\mathcal{I}} \subseteq K$ as the set of nonnegative combinations of $\{\Lambda_i\}_{i \in \mathcal{I}}$. Note that $F_{\mathcal{I}}$ fails to be a face of $\mathcal{C}$ if and only if there exist $\mathcal{I}' \subseteq \mathcal{I}$, $\mathcal{I}'' \not \subseteq \mathcal{I}$, $\alpha_i > 0$ with $\sum_{i \in \mathcal{I}'}\alpha_i = 1$, and $\beta_i > 0$ with $\sum_{i \in \mathcal{I}''}\beta_i = 1$ such that
\begin{equation}
\label{eqC}
y \equiv \sum_{i \in \mathcal{I}'}\alpha_iC_i = \sum_{i \in \mathcal{I}''}\beta_iC_i\,.
\end{equation}
On the other hand $F'_{\mathcal{I}}$ fails to be a face of $K$ if and only if there exist $\mathcal{I}' \subseteq \mathcal{I}$, $\mathcal{I}'' \not \subseteq \mathcal{I}$, $\alpha'_i > 0$, and $\beta'_i > 0$ such that 
\begin{equation}
\label{eqK}
x \equiv \sum_{i \in \mathcal{I}'}\alpha'_i\Lambda_i = \sum_{i \in \mathcal{I}''}\beta'_i\Lambda_i\,.
\end{equation}
If (\ref{eqC}) holds, then $T(y) = \sum_{i \in \mathcal{I}'}\alpha_i\Lambda_i = \sum_{i \in \mathcal{I}''}\beta_i\Lambda_i$ and so, setting $x = T(y)$, $\alpha'_i = \alpha_i$ and $\beta'_i = \beta_i$, (\ref{eqK}) holds and $F'_{\mathcal{I}}$ fails to be a face of $K$. Conversely, suppose (\ref{eqK}) holds. Since $q \not \in \mathrm{Im}\,\Gamma$, there exists some $p \in \mathrm{ker}(\Gamma^{\mathrm{T}})$ such that $p^{\mathrm{T}}q \ne 0$ (choose $p = (2,1,1,1)^{\mathrm{T}}$ as above for example). Writing $x = rT(y)$ and multiplying (\ref{eqK}) from the left by $p^{\mathrm{T}}$ gives
\begin{equation}
\label{eqr}
r = \sum_{i \in \mathcal{I}'}\alpha'_i = \sum_{i \in \mathcal{I}''}\beta'_i\,.
\end{equation}
Noting that $r \ne 0$ and $\mathrm{ker}\,\Gamma = \{0\}$, (\ref{eqr}) and (\ref{eqK}) together imply:
\[
y = \frac{1}{\left(\sum_{i \in \mathcal{I}'}\alpha'_i\right)}\sum_{i \in \mathcal{I}'}\alpha'_iC_i = \frac{1}{\left(\sum_{i \in \mathcal{I}''}\beta'_i\right)}\sum_{i \in \mathcal{I}''}\beta'_iC_i\,.
\]
Defining $\alpha_i = \alpha'_i/(\sum_{i \in \mathcal{I}'} \alpha'_i)$, $\beta_i = \beta'_i/(\sum_{i \in \mathcal{I}''} \beta'_i)$, we see that (\ref{eqC}) holds and $F_{\mathcal{I}}$ fails to be a face of $\mathcal{C}$. More general applications of this argument, and examples of the use of such cones in the study of dynamical systems, can be found in \cite{banajimierczynski}. \\

{\bf The nondegeneracy condition is fulfilled.} That $\mathrm{Im}\,A$ does not lie in $\mathrm{span}\,F$ for any nontrivial face $F$ can be confirmed theoretically, or checked by demonstrating for each three dimensional face $F$ some vector $z$ such that $Az \not \in \mathrm{span}\,F$. This is left to the reader. \\

{\bf $AB$ is $K$-quasipositive.} Finally, that $AB$ is $K$-quasipositive for all $B \in \mathcal{Q}_0(-A^\mathrm{T})$ can easily be checked. Each $B\in \mathcal{Q}_0(-A^\mathrm{T})$ has the form
\[
B= \left(\begin{array}{rrrr}a&-b&0&-c\\0&d&-e&0\\0&0&f&-g\end{array}\right)
\]
where $a,b,c,d,e,f,g \geq 0$. Defining the nonnegative matrix
\[
Q = \left(\begin{smallmatrix} \begin{smallmatrix}a+b+c\\+d+g\end{smallmatrix}&g&d&0&a+b+c&0&0&0\\f&\begin{smallmatrix}a+b+d\\+e+f\end{smallmatrix}&0&e+d&0&a+b&0&0\\e&0&\begin{smallmatrix}a+c+e\\+f+g\end{smallmatrix}&f+g&0&0&a+c&0\\\\0&0&0&a&0&0&0&a\\\\0&0&0&0&0&0&0&0\\\\0&c&0&0&f+g&\begin{smallmatrix}c+e\\+f+g\end{smallmatrix}&0&e\\0&0&b&0&d+e&0&\begin{smallmatrix}b+d\\+e+f\end{smallmatrix}&f\\0&0&0&b+c&0&d&g&\begin{smallmatrix}b+c\\+d+g\end{smallmatrix}\end{smallmatrix}\right)
\]
we can confirm that
\[
AB\Lambda + (a+b+c+d+e+f+g)\Lambda = \Lambda Q.
\]
In other words $AB$ is $K$-quasipositive.


\bibliographystyle{unsrt}

\end{document}